\documentclass[12pt,a4paper]{amsart}

\usepackage{euscript,amsfonts,amssymb,amsmath,amscd,amsthm,enumerate,hyperref}
\usepackage{comment}
\usepackage{mathrsfs}
\usepackage{graphicx}
\usepackage{color}
\usepackage[margin=1.1in]{geometry}
\usepackage{bbm}
\usepackage{upgreek}
\usepackage{enumerate}

\newcommand{\ev}{\mathbb{E}}

\newcommand{\R}{\mathbb{R}}

\newcommand{\dd}{\mathrm{d}}
\newcommand{\pp}{\partial}

\newcommand{\veps}{\varepsilon}

\newcommand{\bP}{\mathbb{P}}

\usepackage{amsthm}

\theoremstyle{plain}
\newtheorem{thm}{Theorem}[section]
\newtheorem{prop}[thm]{Proposition}
\newtheorem{coro}[thm]{Corollary}
\newtheorem{lem}[thm]{Lemma}

\theoremstyle{definition}

\title[Supercooled Stefan problem]{Zero kinetic undercooling limit in the supercooled Stefan problem}

\author{Graeme Baker}
\address{Program in Applied \& Computational Mathematics, Princeton University, Princeton, NJ 08544, USA}
\email{graemeb@princeton.edu}

\author{Mykhaylo Shkolnikov}
\address{ORFE Department, Bendheim Center for Finance, and Program in Applied \& Computational Mathematics, Princeton University, Princeton, NJ 08544, USA}
\email{mshkolni@gmail.com}
\thanks{G.~Baker is partially supported by an NSERC PGS-D scholarship and a Princeton SEAS innovation research grant. M.~Shkolnikov is partially supported by the NSF grant DMS-1811723 and a Princeton SEAS innovation research grant.}

\begin{document}

\begin{abstract}
We study the solutions of the one-phase supercooled Stefan problem with kinetic undercooling, which describes the freezing of a supercooled liquid, in one spatial dimension. Assuming that the initial temperature lies between the equilibrium freezing point and the characteristic invariant temperature throughout the liquid our main theorem shows that, as the kinetic undercooling parameter tends to zero, the free boundary converges to the (possibly irregular) free boundary in the supercooled Stefan problem without kinetic undercooling, whose uniqueness has been recently established in \cite{delarue2019global}, \cite{LS}. The key tools in the proof are a Feynman-Kac formula, which expresses the free boundary in the problem with kinetic undercooling through a local time of a reflected process, and a  resulting comparison principle for the free boundaries with different kinetic undercooling parameters. 
\end{abstract}

\maketitle

%%%%%%%%%%%%%
\section{Introduction}\label{intro}
%%%%%%%%%%%%%
%history of Stefan
Free boundary problems for the heat equation which describe the freezing of liquids have been a subject of research since the time of \textsc{Lam\'e} and \textsc{Clapeyron} \cite{LC}. Such problems are named after \textsc{Stefan}, as he studied them systematically in the series of papers \cite{Stefan1,Stefan2,Stefan3,Stefan4} and, in particular, solved the problem for a liquid in the half-space $\{x_1>0\}$ with a constant temperature below the equilibrium freezing point maintained at the boundary $\{x_1=0\}$ and instantaneous freezing of the liquid at the equilibrium freezing point. The lecture of \textsc{Brillouin} in 1929 reinitiated the investigation of Stefan problems after a period of dormancy, with a focus on existence, uniqueness, and numerical approximation (see \cite{Bri}). In \cite{Kam}, \textsc{Kamenomostskaja} established existence and uniqueness and gave an explicit difference scheme for bounded measurable generalized solutions in any dimension and for any number of phases.

\medskip

%history of supercooled
Less is known for the \emph{supercooled} Stefan problem, in which the initial temperature of the liquid lies below its equilibrium freezing point. This is due to the fact that already the one-phase problem in one spatial dimension can exhibit a finite time \emph{blow-up} of the freezing rate,  as first noted in \cite{Sher}. One workaround, considered in \cite{Vis,DHOX,HX,FPHO2,Wei} and referred to as kinetic undercooling, is to regularize the problem and prevent blow-up through a modification of the boundary condition. On the other hand, the probabilistic reformulation of the original supercooled Stefan problem in \cite{delarue2019global} (see also \cite{HLS,ledger2018mercy,LS,nadtochiy2018mean,NS}) allows to examine the global solutions even in the presence of blow-ups. These probabilistic solutions are known to exist (\cite[Theorem 3.2]{ledger2018mercy},  \cite[Theorem 2.3]{nadtochiy2018mean}, \cite[Theorem 2.4]{NS}) and be unique (\cite[Theorem 1.4]{delarue2019global}, \cite[Theorem 2.2]{LS}) under mild assumptions on the initial temperature distribution, and the associated free boundary to transition between continuous differentiability, H\"older continuity, and discontinuity (\cite[Theorem 1.1]{delarue2019global}). 

\medskip

%justify global solution and connect with regularized version
The goal of this paper is to demonstrate that the solutions of the supercooled Stefan problem with kinetic undercooling are intimately connected to the probabilistic solutions of the supercooled Stefan problem without kinetic undercooling. More specifically, we prove in Theorem \ref{thm:main} that the free boundaries of the regularized problems converge uniformly on compact sets to the probabilistic free boundary in the problem without regularization, as the regularization parameter tends to zero. The link between the problems is made clear by a new probabilistic reformulation of the regularized problems. In light of our findings, the probabilistic solution is justifiably deemed to be an appropriate global notion of solution for the supercooled Stefan problem without kinetic undercooling.

\medskip

The one-dimensional one-phase supercooled Stefan problem is the free boundary problem
\begin{equation}\label{stefan}
\begin{cases}
\; \pp_tu=\frac{1}{2}\pp_{xx}u\;\;\text{on}\;\;\Gamma:=\{(t,x)\in[0,\infty)^2:x\ge\Lambda(t)\}, \\
\; u(0,\cdot)=f\;\;\text{and}\;\;u(t,\Lambda(t))=0,\;t\ge0, \\ 
\;\alpha\Lambda'(t)=\pp_x u(t,\Lambda(t)),\;t\ge0\;\;\text{and}\;\;\Lambda(0)=0,
\end{cases}
\end{equation}
with non-negative initial data $f$ and a constant $\alpha>0$ representing the density of latent heat. %(note that \cite{delarue2019global} uses $2/\alpha$ in place of $\alpha$).
Hereby, $(\Lambda(t),\infty)$ captures the domain occupied by the supercooled liquid at time $t\ge0$, and $u(t,x)$ gives the number of degrees below the equilibrium freezing point at time $t\ge0$ and position $x\ge\Lambda(t)$. The heat capacity of the liquid is set to 2 and its thermal conductivity to 1, resulting in the (probabilist's) standard heat equation $\pp_tu=\frac{1}{2}\pp_{xx}u$. 

\medskip

The regularized problem with kinetic undercooling and regularization parameter $\veps>0$ reads
\begin{equation}\label{stefan_reg}
\begin{cases}
\;\pp_tu_\veps=\frac{1}{2}\pp_{xx}u_\veps\;\;\text{on}\;\;\Gamma_\veps:=\{(t,x)\in[0,\infty)^2:x\ge\Lambda_\veps(t)\}, \\
\;u_\veps(0,\cdot)=f_\veps\;\;\text{and}\;\;u_\veps(t,\Lambda_\veps(t))=\veps\Lambda'_\veps(t),\;t\ge0, \\ 
\;(\alpha-2u_\veps(t,\Lambda_\veps(t)))\,\Lambda'_\veps(t)=\pp_x u_\veps(t,\Lambda_\veps(t)),\;t\ge0\;\;\text{and}\;\;
\Lambda_\veps(0)=0,
\end{cases}
\end{equation}
where $f_\veps := f*\rho_\veps$, $\rho_\veps(\cdot):=\rho(\cdot/\veps)/\veps$, and $\rho\in C_c^\infty((0,\infty))$ is a non-negative smoothing kernel that integrates to one.
%The solid is kept at constant equilibrium temperature 
Let us explain the origins of these equations. The kinetic undercooling equation assumes that the liquid-to-solid phase transition is driven by undercooling, leading to $\Lambda'_\veps(t)=\gamma_\veps(u_\veps(t,\Lambda_\veps(t)))$ and after linearization of the kinetic function $\gamma_\veps$ around $u_\veps=0$ to $\Lambda'_\veps(t)=\frac{1}{\veps}u_\veps(t,\Lambda_\veps(t))$.~The Stefan condition states $\pp_x u_\veps(t,\Lambda_\veps(t))=L(u_\veps(t,\Lambda_\veps(t)))\,\Lambda'_\veps(t)$ (see, e.g., \cite[equation (1.2)]{visintin1998introduction}) and following \cite[equation (2.23)]{glicksman2010principles} we use a linear approximation for the density of latent heat $L$ to obtain $\pp_x u_\veps(t,\Lambda_\veps(t))=(\alpha-2u_\veps(t,\Lambda_\veps(t)))\,\Lambda'_\veps(t)$. We note that the problem with kinetic undercooling \eqref{stefan_reg} is likely only physically meaningful when $\alpha-2u_\veps>0$ since the density of latent heat becomes zero at the characteristic invariant temperature $-u_\veps=-\frac{\alpha}{2}$, beyond which (in the so-called hypercooled regime) spontaneous nucleation and liquid-glass transition have been observed in physical experiments (see \cite[Subsection 17.3.2]{glicksman2010principles}). Therefore, we take $\|f\|_\infty<\frac{\alpha}{2}$ throughout. 

In a nutshell, the absence of blow-ups in \eqref{stefan_reg} is due to $\Lambda_\veps'(t)=\frac{1}{\veps}u_\veps(t,\Lambda_\veps(t))$, $t\ge0$, where the non-negative $u_\veps$ cannot exceed $\|f_\veps\|_\infty\le\|f\|_\infty<\frac{\alpha}{2}$ by a maximum principle. Thus, $\Lambda'_\veps$ is controlled a priori by $\|f_\veps\|_\infty/\veps$ for all times, ruling out blow-ups. We refer to the proof of Proposition \ref{PDEsol} below for more details.
%As $\veps\downarrow 0$, we will see that the regularized free boundary $\Lambda_\veps$ converges in the Skorokhod M1 topology to the boundary from \cite{delarue2019global}; this is stated precisely  below as Theorem \ref{thm:main}.

\medskip

Next, we describe the probabilistic reformulation of \eqref{stefan} from \cite{delarue2019global} (see also \cite{HLS}, \cite{LS}) that allows to make sense of \eqref{stefan} globally, despite the possible blow-ups of $\Lambda'$. Suppose $\int_0^\infty f(x)\,\mathrm{d}x=1$ and let $X_{0-}\geq 0$ be a random variable with density $f$. For a standard Brownian motion $B$ independent of $X_{0-}$, consider the problem of finding a non-decreasing right-continuous function $\Lambda:\,[0,\infty)\to[0,\infty)$ such that
\begin{align}
\label{limprob}
\begin{cases}
\;\alpha\Lambda(t)=2\bP(\tau\leq t),\;t\geq0,\;\;\text{where} \\
\;\tau:=\inf\{t\geq0:\; X_t\leq 0\} \;\;\text{and}\;\; X_t:=X_{0-}+B_t-\Lambda(t),\;t\geq0.
\end{cases}
\end{align}
Writing $p(t,x)\,\mathrm{d}x$, $t\ge0$ for the distributions of $X_t\,\mathbf{1}_{\{\tau>t\}}$, $t\ge 0$ restricted to $(0,\infty)$, the function $u(t,x):=p(t,x-\Lambda(t))$, $(t,x)\in\Gamma$ combines with
the solution $\Lambda$ of \eqref{limprob} to a global solution of \eqref{stefan}, as explained in the introduction of  \cite{delarue2019global}. Under our standing assumption $\|f\|_\infty<\frac{\alpha}{2}$ the solution $\Lambda$ of \eqref{limprob} is unique and continuous. 

\begin{prop}[\cite{LS}, Theorem 2.2 and last paragraph of Section 2.1] \label{thm:unique} 
If the density $f$ of $X_{0-}$ obeys $\|f\|_\infty<\frac{\alpha}{2}$, then the solution $\Lambda$ of \eqref{limprob} is unique and continuous. 
\end{prop}

We can now state our main result. 

\begin{thm}\label{thm:main}
Let a non-negative bounded initial data $f$ be given with $\int_0^\infty f(x)\,\mathrm{d}x=1$ and $\|f\|_\infty<\frac{\alpha}{2}$. Then, for each $\veps>0$, the problem \eqref{stefan_reg} admits a unique free boundary $\Lambda_\veps$ with $\|\Lambda'_\veps\|_\infty\le\|f_\veps\|_\infty/\veps$. The family $\{\Lambda_\veps\}_{\veps>0}$ increases pointwise to the unique solution $\Lambda$ of \eqref{limprob} as $\veps\downarrow0$ and, thus, converges uniformly on compact sets to the latter. 
\end{thm}

The proof of Theorem \ref{thm:main} serves as a roadmap to the rest of the paper.

\medskip

\noindent\textbf{Proof of Theorem \ref{thm:main}.}
In Section \ref{sec:aux} we introduce an auxiliary fixed boundary PDE \eqref{ppde} and cast \eqref{stefan_reg} as a solution to \eqref{ppde} which satisfies a fixed point problem.
The key tool for working with the boundaries in \eqref{stefan_reg} and \eqref{ppde} is a Feynman-Kac formula developed in Section \ref{sec:FK}. The existence and uniqueness of the free boundaries $\{\Lambda_\veps\}_{\veps>0}$ is shown via Banach's fixed point theorem in Proposition \ref{PDEsol}. In Proposition \ref{prop:mono} we prove the pointwise monotonicity of $\{\Lambda_\veps\}_{\veps>0}$ in $\veps$. In Proposition \ref{prop:limprob} we check that the right-continuous modification of the limit $\widetilde{\Lambda}$ solves the problem \eqref{limprob}, and in Proposition \ref{prop:unique} we verify that $\widetilde{\Lambda}$ agrees with its unique solution $\Lambda$. \qed

\section{Auxiliary Problem}\label{sec:aux}
In this section we set up the auxiliary fixed boundary problem \eqref{ppde} which is central to our investigations. For $T\in(0,\infty)$, write $\widetilde{W}^1_\infty([0,T])$ for the Banach space of bounded Lipschitz functions from $[0,T]$ to $\R$ vanishing at the left endpoint $0$, with the Lipschitz constant as the norm. Now, for $\Lambda_\veps\in\widetilde{W}^1_\infty([0,T])$, consider the problem   
\begin{equation}\label{ppde}
\begin{cases}
\; \pp_tp_\veps=\frac{1}{2}\pp_{xx}p_\veps+\Lambda'_\veps\,\pp_x p_\veps\;\;\text{on}\;\;[0,T]\times[0,\infty), \\
\; p_\veps(0,\cdot)=f_\veps\;\;\text{and}\;\; 
\pp_x p_\veps(\cdot,0)=\left(\frac{\alpha}{\veps}-2\Lambda_\veps'\right)p_\veps(\cdot,0). 
\end{cases}
\end{equation}
For any non-negative bounded $f$ with $\int_0^\infty f(x)\,\mathrm{d}x=1$, there exists a unique solution $p_\veps$ in the Sobolev space $W^{1,2}_2([0,T]\times[0,\infty))$ (see \cite[Chapter IV, Section 9]{lady1968linear}: our problem corresponds to (5.4) therein, $\Lambda_\veps'$ is bounded, and $f'_\veps(0)=f_\veps(0)=0$).
Further, by \cite[Chapter II, Lemma 3.3]{lady1968linear} the solution $p_\veps$ admits a (H\"older) continuous version, allowing us to define the following operator: 
\begin{align}
F_\veps:\;\widetilde{W}^1_\infty([0,T])\to\widetilde{W}^1_\infty([0,T]),\quad
\Lambda_\veps\mapsto\frac{1}{\veps}\int_0^\cdot p_\veps(s,0)\,\dd s. \label{Fdef}
\end{align}

\smallskip

Suppose that $\Lambda_\veps\in \widetilde{W}^1_\infty([0,T])$ is a fixed point of $F_\veps$ and set $u_\veps(t,x)=p_\veps(t,x-\Lambda_\veps(t))$, $x\ge\Lambda_\veps(t)$, $t\in[0,T]$. Then, $\partial_t u_\veps=\frac{1}{2}\partial_{xx}u_\veps$ a.e., $u_\veps(0,\cdot)=f_\veps$, and $\partial_xu_\veps(t,\Lambda_\veps(t))=\left(\frac{\alpha}{\veps}-2\Lambda_\veps'(t)\right)u_\veps(t,\Lambda_\veps(t))$, $t\in[0,T]$. In addition, differentiating \eqref{Fdef} we get $\Lambda'_\veps(t)=\frac{1}{\veps}\,p_\veps(t,0)=\frac{1}{\veps}\,u_\veps(t,\Lambda_\veps(t))$, $t\in[0,T]$, and together with the Robin boundary condition:
\begin{equation}
\Lambda_\veps'(t)= \frac{1}{\alpha}\big(\partial_xu_\veps(t,\Lambda_\veps(t))+2\Lambda_\veps'(t)\,u_\veps(t,\Lambda_\veps(t))\big),\;t\in[0,T].
\end{equation}
Thus, any fixed point of $F_\veps$ in $\widetilde{W}^1_\infty([0,T])$ leads to a solution of \eqref{stefan_reg} on $[0,T]$. In Section \ref{sec:PDE}, we show that $F_\veps$ exhibits a unique fixed point $\Lambda_\veps$ with $\|\Lambda'_\veps\|_\infty\le\|f_\veps\|_\infty/\veps$, but first we develop some useful tools in the next section.

%\section{Probabilistic Formulation for the Regularized Problem}\label{sec:FK}
%%%%%%%%%%%%%%%%%%%
\section{Feynman-Kac Formula}\label{sec:FK}
%%%%%%%%%%%%%%%%%%%

In this section, we develop the main tool for our purposes, a Feynman-Kac formula for the problem \eqref{ppde}. The latter should be viewed in analogy to the probabilistic formulation \eqref{limprob} for the problem without regularization. Instead of the absorbed process $X_t\,\mathbf{1}_{\{\tau>t\}}$, $t\ge0$ therein, we employ the \textit{reflection} $X^\veps$ of the process $X_0^\veps+B_t-\Lambda_\veps(t)$, $t\in[0,T]$ at zero, initialized according to the density $f_\veps$. The role of $\mathbf{1}_{\{\tau\leq t\}}$ is now played by $1-e^{-L_t^{\veps,0}/\veps}$, where $L^{\veps,0}$ is the local time of $X^\veps$ at $0$. With $\tau^\veps:=\inf\{t\geq0:\,X^\veps_t=0\}$, we have $1-e^{-L_t^{\veps,0}/\veps}=\mathbf{1}_{\{\tau^\veps\le t \}}\,(1-e^{-L_t^{\veps,0}/\veps})$, and so we expect, heuristically, that $1-e^{-L_t^{\veps,0}/\veps}\approx\mathbf{1}_{\{\tau^\veps\leq t\}}$ as $\veps\to 0$. The rigorous analysis builds on the Skorokhod lemma (see, e.g., \cite[Chapter 3, Lemma 6.14]{karatzas1998brownian}).  

\begin{lem}\label{Skoro}
For any $T\in(0,\infty)$ and continuous $y:\,[0,T]\to\R$, there exists a unique continuous  $\ell:\,[0,T]\to\R$ such that
\begin{enumerate}[(i)]
\item $\Psi(y)(t):=y(t)+\ell(t)\geq0$, $t\in[0,T]$,
\item $\ell(0)=0$ and $\ell$ is non-decreasing, 
\item $\int_0^T \mathbf{1}_{\{\Psi(y)(s)>0\}}\,\dd\ell(s)=0$.
\end{enumerate}
Moreover, this function admits the explicit formula
\begin{align}\label{SkorFormula}
\ell(t)=0\vee \big(\!-\min_{0\leq s \leq t}\,y(s)\big),\;t\in[0,T].
\end{align}
The map $\Psi:\,C([0,T])\to C([0,T])$ is referred to as the  \emph{Skorokhod map}.
\end{lem}

We apply the Skorokhod map to the paths of the process $X_0^\veps+B_t-\Lambda_\veps(t)$, $t\in[0,T]$.

\begin{prop}\label{prop:refl}
For non-negative $f$ with $\int_0^\infty f(x)\,\mathrm{d}x=1$, let $X^\veps=\Psi(X^\veps_0+B-\Lambda_\veps)$, where $X^\veps_0$ is a random variable with density $f_\veps$ and $\Lambda_\veps\in\widetilde{W}^1_\infty([0,T])$. Then, $X^\veps$ satisfies the stochastic differential equation 
\begin{align}
\dd X^\veps_t=\dd B_t -\Lambda_\veps'(t)\,\dd t+\dd L^{\veps,0}_t
\end{align}
on $[0,T]$, with the local time $L^{\veps,0}$ of $X^\veps$ at $0$. 
\end{prop}

\noindent\textbf{Proof.} The result for $\Lambda_\veps\equiv 0$ is well-known (see, e.g., \cite[Chapter 3, display (6.33) and the subsequent sentence]{karatzas1998brownian}). The general case is readily obtained via Girsanov's theorem, noticing that $\Lambda'_\veps$ is bounded by assumption. \qed

\medskip

%Now, we recall the notion of local time and the Skorokhod lemma for reflecting processes of the form \eqref{eq:reflSDE}.
%Our main reference is Chapter VI of \cite{revuz2013continuous}.

%\begin{defn}[Theorem VI.1.2 in \cite{revuz2013continuous}]
%	For any continuous semimartingale $X$ and any $x\in\R$ there exists an non-decreasing continuous process $L^x$ called the \emph{local time of $X$ at $x$} such that
%	\begin{align}
%	(X_t-x)^+=(X_0-x)^+ +\int_0^t 1_{\{X_s\geq x\}}\dd X_s+\frac{1}{2}L_t^x.
%	\end{align}
%\end{defn}

We collect some properties of local time for future use. 

\begin{prop}\label{prop:lt}
For a continuous semimartingale $X$ with local time $L$ it holds:
\begin{enumerate}[(a)]
\item \label{prop:lt:cont} There exists a modification of $L$ such that $(t,x)\to L_t^x$ is a.s.~continuous in $t$ and right-continuous with left limits in $x$. This is the version we consider throughout. 
\item A.s., for every $t$ and $x$, 
\begin{align}
L_t^x=\lim_{\delta\downarrow 0}\,\frac{1}{\delta}\int_0^t \mathbf{1}_{[x,x+\delta)}(X_s)\,\dd \langle X,X\rangle_s.
\end{align}
\item \label{prop:lt:supp} For each $x$, the support of the measure $\dd L^x$ is contained in $\{t\ge0:\,X_t=x\}$ a.s.
%equal to this when X is BM. What about BM with drift? Where does drift have to be?
\item \label{prop:lt:occ} (Occupation time formula) A.s., for any non-negative Borel $g$,
\begin{align}
\int_0^t g(s,X_s)\,\dd \langle X,X\rangle_s = \int_{-\infty}^\infty\int_0^t g(s,x)\,\dd L_s^x\,\dd x.
\end{align}
\noindent For $X^\veps=\Psi(X^\veps_0+B-\Lambda_\veps)$ as in Proposition \ref{prop:refl}, we also have: \smallskip
\item \label{prop:lt:holder} The local time $L^\veps$ may be chosen so that, a.s., upon replacing $L^{\veps,0}$ by $2L^{\veps,0}$ the functions $x\mapsto L^{\veps,x}_t$ become $\gamma$-H\"older continuous for all $\gamma<1/2$ uniformly over compact intervals in $t$.
\item \label{prop:lt:sup} $\ev\big[(\sup_{x\in\R} L_T^{\veps,x})^p\big]<\infty$, $p\in(0,\infty)$. 
\end{enumerate}
\end{prop}

\noindent\textbf{Proof.} Parts (\ref{prop:lt:cont})-(\ref{prop:lt:occ}) can be found in Theorem 1.7, Corollary 1.9, Proposition 1.3, and Exercise 1.15 of \cite[Chapter VI]{revuz2013continuous}, respectively. For $\Lambda_\veps\equiv0$,  parts (\ref{prop:lt:holder}), (\ref{prop:lt:sup}) follow from \cite[Chapter 3, Proposition 6.16]{karatzas1998brownian} and \cite[Chapter VI, Corollary 1.8]{revuz2013continuous}, \cite[Chapter X, Exercise 1.14]{revuz2013continuous}, respectively. For  $\Lambda_\veps\in \widetilde{W}^1_\infty([0,T])$, it suffices to employ Girsanov's theorem. \qed. 

\medskip

We are now in a position to state our Feynman-Kac formula.

\begin{prop}\label{prop:FK}
For a non-negative bounded $f$ with $\int_0^\infty f(x)\,\mathrm{d}x=1$, a random variable $X_0^\veps$ with density $f_\veps$, a function $\Lambda_\veps\in\widetilde{W}^1_\infty([0,T])$ and $X^\veps:=\Psi(X_0^\veps+B-\Lambda_\veps)$, write $\widetilde{p}_\veps(t,\cdot)$, $t\in[0,T]$ for the densities of the sub-probability measures given by 
\begin{equation}
\mu_t^\veps([a,b))= \ev\big[\mathbf{1}_{\{X_t^\veps\in [a,b)\}}\,e^{-\alpha L_t^{\veps,0}/\veps}\big],\;0\le a<b<\infty,\; t\in[0,T]. 
\end{equation}
Then, for all $t\in[0,T]$, one has $\widetilde{p}_\veps(t,\cdot)=p_\veps(t,\cdot)$ a.e., where $p_\veps$ is the solution of \eqref{ppde} in $W^{1,2}_2([0,T]\times[0,\infty))$.
\end{prop}

\noindent\textbf{Proof.} Fix a $t\in(0,T]$ and a $\xi_\veps\in C^1_c((0,\infty))$, and let $\eta_\veps\in W^{1,2}_2([0,t]\times [0,\infty))$ be the unique solution of
\begin{align}
\begin{cases}\label{eq:eta}
\pp_s\eta_\veps+\frac{1}{2}\pp_{xx}\eta_\veps-\Lambda'_\veps\,\pp_x\eta_\veps=0\;\;\text{on}\;\;[0,t]\times[0,\infty), \\
\eta_\veps(t,\cdot)=\xi_\veps\;\;\text{and}\;\;\alpha\eta_\veps(\cdot,0)=\veps\pp_x\eta_\veps(\cdot,0)
\end{cases}
\end{align}
(see \cite[Chapter IV, Section 9]{lady1968linear}). Since $e^{-\alpha x/\veps}\eta_\veps\in W^{1,2}_2([0,t]\times [0,\infty))$ and $\partial_x(e^{-\alpha x/\veps}\eta_\veps)(\cdot,0)\equiv 0$, we can apply It\^o's formula to $e^{-\alpha \Psi(X_0^\veps+B)_s/\veps}\eta_\veps(s,\Psi(X_0^\veps+B)_s)$ (see \cite[Section 2.10, proof of Theorem 1]{krylov2008controlled}), which with the product rule and Girsanov's theorem yields
\begin{equation}
\begin{split}\label{prod_Ito}
&\,e^{-\alpha L_t^{\veps,0}/\veps}\eta_\veps(t,X_t^\veps)-\eta_\veps(0,X_0^\veps)\\
&=\int_0^te^{-\alpha L_s^{\veps,0}/\veps} \Big(\pp_s \eta_\veps(s,X_s^\veps)+\frac{1}{2}\pp_{xx}\eta_\veps(s,X_s^\veps)-\Lambda'_\veps(s)\,\pp_x \eta_\veps(s,X_s^\veps)\Big)\,\dd s\\
&\quad+\int_0^t e^{-\alpha L_s^{\veps,0}/\veps}\pp_x \eta_\veps(s,X_s^\veps)\,\dd B_s + \int_0^t e^{-\alpha L_s^{\veps,0}/\veps}\Big(\pp_x \eta_\veps(s,X_s^\veps)-\frac{\alpha }{\veps}\eta_\veps(s,X_s^\veps)\Big)\,\dd L_s^{\veps,0} \\
&=\int_0^t e^{-\alpha L_s^{\veps,0}/\veps}\pp_x \eta_\veps(s,X_s^\veps)\,\dd B_s,
\end{split}
\end{equation}
where the second equality is due to \eqref{eq:eta} and Proposition \ref{prop:lt}(\ref{prop:lt:supp}). 

\medskip
	
We infer $\ev\big[\int_0^t e^{-2\alpha L_s^{\veps,0}/\veps}\pp_x \eta_\veps(s,X_s^\veps)^2\,\dd s\big]<\infty$ using Girsanov's theorem and twice the Cauchy-Schwarz inequality, and arguing $\ev\big[\int_0^t \pp_x \eta_\veps(s,\Psi(X_0^\veps+B)_s)^4\,\dd s\big]<\infty$ via the boundedness of the densities of $\Psi(X_0^\veps+B)_s$, $s\in[0,t]$ by $\|f_\veps\|_\infty$ (note that the transition kernel of $\Psi(X_0^\veps+B)$ %, given in \eqref{trans},
integrates to $1$ in the backward variable) and $\pp_x \eta_\veps\in L^4([0,t]\times[0,\infty))$ (see \cite[Chapter II, Lemma 3.3]{lady1968linear}). Taking the expectation in \eqref{prod_Ito} and recalling $\eta_\veps(t,\cdot)=\xi_\veps$ as well as the definition of $\widetilde{p}_\veps(t,\cdot)$, we arrive at  
\begin{align}
\label{eq:tildep}
\int_0^\infty \xi_\veps(x)\,\widetilde{p}_\veps(t,x)\,\dd x=\int_0^\infty \eta_\veps(0,x)\,f_\veps(x)\,\dd x.
\end{align}

\smallskip

On the other hand, the PDE in \eqref{ppde} and repeated integration by parts give
\begin{equation}
\begin{split}
0&=\int_0^t\int_{0}^\infty \Big(\pp_s p_\veps(s,x)-\frac{1}{2}\pp_{xx}p_\veps(s,x)
-\Lambda'_\veps(s)\,\pp_x p_\veps(s,x)\Big)\,g(s,x)\,\dd x\,\dd s\\
&=\int_{0}^\infty g(t,x)\,p_\veps(t,x)\,\dd x-\int_0^\infty g(0,x)\,f_\veps(x)\,\dd x\\
&\quad-\int_0^t \int_{0}^\infty \Big(\pp_s g(s,x)+\frac{1}{2}\pp_{xx}g(s,x)-\Lambda'_\veps(s)\,
\pp_x g(s,x)\Big)\,p_\veps(s,x)\,\dd x\,\dd s\\
&\quad+\int_0^t g(s,0)\,\Big(\frac{1}{2}\pp_x p_\veps(s,0)+\Lambda_\veps'(s)\,p_\veps(s,0)\Big)-\frac{1}{2}\pp_xg(s,0)\,p_\veps(s,0)\,\dd s,
\end{split}
\end{equation}
for all $g\in W^{1,2}_2([0,t]\times[0,\infty))$. Strictly speaking, to obtain the second equality we approximate $p_\veps$ and $g$ in the $W^{1,2}_2([0,t]\times[0,\infty))$-norm by $C_c^\infty([0,t]\times[0,\infty))$-functions, perform the integrations by parts on the latter, and pass to the limit relying on the continuity of the evaluation maps in \cite[Chapter II, Lemma 3.4]{lady1968linear}. Thanks to the boundary condition in \eqref{ppde},
\begin{equation}\label{duality}
\begin{split}
&\,\int_{0}^\infty g(t,x)\,p_\veps(t,x)\,\dd x-\int_0^\infty g(0,x)\,f_\veps(x)\,\dd x\\
&=\int_0^t \int_{0}^\infty \Big(\pp_s g(s,x)+\frac{1}{2}\pp_{xx}g(s,x)-\Lambda'_\veps(s)\,
\pp_x g(s,x)\Big)\,p_\veps(s,x)\,\dd x\,\dd s\\
&\quad+\int_0^t \Big(\frac{1}{2}\pp_x g(s,0)-\frac{\alpha}{2\veps}g(s,0)\Big)\,p_\veps(s,0)\,\dd s.
\end{split}
\end{equation}
Plugging in $\eta_\veps$ from above for $g$ yields
\begin{align}
\int_0^\infty \xi_\veps(x)\,p_\veps(t,x)\,\dd x=\int_0^\infty \eta_\veps(0,x)\,f_\veps(x)\,\dd x.
\end{align}
We conclude by recalling \eqref{eq:tildep} and the arbitrariness of  $\xi_\veps\in C^1_c((0,\infty))$. \qed

\begin{coro}\label{sformula}
In the setting of the previous proposition one has the representation
\begin{align}\label{eq:sformula}
F_\veps(\Lambda_\veps)(t):=\frac{1}{\veps} \int_0^t p_\veps(s,0)\,\dd s=\frac{2}{\alpha}  \big(1-\mathbb{E}\big[ e^{-\alpha L_t^{\veps,0}/\veps}\big]\big),\;t\in[0,T].
\end{align}
\end{coro}

\noindent\textbf{Proof.} Using the uniform continuity of $p_\veps$ on $[0,T]\times[0,1]$ we obtain
\begin{align}
F_\veps(\Lambda_\veps)(t) &=\lim_{\delta\downarrow 0}\frac{1}{\veps\delta} \int_0^t \int_0^\delta p_\veps(s,x)\,\dd x\,\dd s,\;t\in[0,T].
\end{align}
Proposition \ref{prop:FK} allows us to write
\begin{equation}
\begin{split}
	\lim_{\delta\downarrow 0}\frac{1}{\veps\delta} \int_0^t \int_0^\delta p_\veps(s,x)\,\dd x\,\dd s&=\lim_{\delta\downarrow 0}\frac{1}{\veps\delta} \int_0^t \int_0^\delta \widetilde{p}_\veps(s,x)\,\dd x\,\dd s\\
	&=\lim_{\delta\downarrow 0}\frac{1}{\veps\delta}\int_0^t\mathbb{E}\big[\mathbf{1}_{\{X_s^\veps\in [0,\delta)\}}\,e^{-\alpha L_s^{\veps,0}/\veps}\big]\,\dd s\\
	&=\lim_{\delta\downarrow 0}\frac{1}{\veps\delta}\mathbb{E}\bigg[\int_0^t \mathbf{1}_{\{X_s^\veps\in [0,\delta)\}}\,e^{-\alpha L_s^{\veps,0}/\veps}\,\dd s\bigg],
\end{split}
\end{equation}
with the last equality being due to Fubini's theorem. The occupation time formula (Proposition \ref{prop:lt}(\ref{prop:lt:occ})) and another instance of Fubini's theorem give
\begin{equation}
\begin{split}
	\lim_{\delta\downarrow 0}\frac{1}{\veps\delta}\mathbb{E}\bigg[\int_0^t \mathbf{1}_{\{X_s^\veps\in [0,\delta)\}}\,e^{-\alpha L_s^{\veps,0}/ \veps}\,\dd s\bigg]
	&=\lim_{\delta\downarrow 0}\frac{1}{\veps\delta}\mathbb{E}\bigg[\int_0^\delta \int_0^te^{-\alpha L_s^{\veps,0}/ \veps}\,\dd L_s^{\veps,x}\,\dd x \bigg]\\
	&=\frac{1}{\veps}\lim_{\delta\downarrow 0}\frac{1}{\delta}\int_0^\delta \mathbb{E}\bigg[\int_0^te^{-\alpha L_s^{\veps,0}/ \veps}\,\dd L_s^{\veps,x}\bigg]\,\dd x.\label{eq:occ}
\end{split}
\end{equation}
Note that $\mathbb{E}\big[\int_0^te^{-\alpha L_s^{\veps,0}/\veps}\,\dd L_s^{\veps,x}\big]$ becomes continuous in $x$ at $0$ once we replace $\mathrm{d}L^{\veps,0}$ by $2\,\mathrm{d}L^{\veps,0}$. Indeed, after the replacement $\int_0^t e^{-\alpha L_s^{\veps,0}/ \veps}\,\dd L_s^{\veps,x}$ becomes continuous in $x$ by Proposition \ref{prop:lt}(\ref{prop:lt:holder}) and the dominated convergence theorem applies since we may bound $e^{-\alpha L_s^{\veps,0}/ \veps}$ by $1$ and $\mathbb{E}[\sup_{x\in\R} L_t^{\veps,x}]<\infty$ by Proposition \ref{prop:lt}(\ref{prop:lt:sup}).

\medskip	
	
We manipulate the end result of \eqref{eq:occ} to obtain the desired representation:
\begin{equation}
	\frac{1}{\veps}\lim_{\delta\downarrow 0}\frac{1}{\delta}\int_0^\delta \mathbb{E}\bigg[\int_0^te^{-\alpha L_s^{\veps,0}/\veps}\,\dd L_s^{\veps,x}\bigg]\,\dd x
	=\frac{2}{\veps} \mathbb{E}\bigg[\int_0^te^{-\alpha L_s^{\veps,0}/\veps}\,\dd L_s^{\veps,0}\bigg]
	=\frac{2}{\alpha} \big(1- \mathbb{E}\big[e^{-\alpha L_t^{\veps,0}/\veps}\big]\big),
\end{equation}
where the final equality holds by explicit integration. \qed

\medskip

We frequently put Corollary \ref{sformula} together with the following comparison principle. 

\begin{prop}\label{prop:burdzy}
For $\veps\ge\widehat{\veps}>0$ and non-negative $f$ with $\int_0^\infty f(x)\,\mathrm{d}x=1$, let $X^\veps=\Psi(X^\veps_0+B-\Lambda_\veps)$ and $\widehat{X}^{\widehat{\veps}}=\Psi(\widehat{X}^{\widehat{\veps}}_0+B-\widehat{\Lambda}_{\widehat{\veps}})$, where $X^\veps_0\ge\widehat{X}^{\widehat{\veps}}_0$ are random variables with densities $f_\veps$, $f_{\widehat{\veps}}$, respectively, and $\Lambda_\veps,\widehat{\Lambda}_{\widehat{\veps}}\in \widetilde{W}^1_\infty([0,T])$ with $\Lambda_\veps\le \widehat{\Lambda}_{\widehat{\veps}}$. Then, $L^{\veps,0}\le \widehat{L}^{\widehat{\veps},0}$ on $[0,T]$, for the local times $L^{\veps,0}$, $\widehat{L}^{\widehat{\veps},0}$ of $X^\veps$, $\widehat{X}^{\widehat{\veps}}$, respectively, at $0$. 
\end{prop}

\noindent\textbf{Proof.} It suffices to combine Proposition \ref{prop:refl} and the formula \eqref{SkorFormula} with $y=X^\veps_0+B-\Lambda_\veps$ and $y=\widehat{X}^{\widehat{\veps}}_0+B-\widehat{\Lambda}_{\widehat{\veps}}$. \qed

%remark on initial conditions
%%%%%%%%%%%%%%%%
\section{Solving the regularized problem}\label{sec:PDE}
%%%%%%%%%%%%%%%%

We now tackle the fixed point problem for $F_\veps$ introduced at the end of Section \ref{sec:aux}. Recall that every fixed point of $F_\veps$ leads to a solution of the regularized problem \eqref{stefan_reg}.
%upgrade to theorem?
\begin{prop}\label{PDEsol}
For non-negative bounded $f$ with $\int_0^\infty f(x)\,\mathrm{d}x=1$ and $\|f\|_\infty\!<\frac{\alpha}{2}$, the map $F_\veps$ possesses a unique fixed point $\Lambda_\veps$ in $\widetilde{W}^{1}_\infty([0,T])$ with $\|\Lambda'_\veps\|_\infty\le\|f_\veps\|_\infty/\veps$. In addition, for all small enough $T>0$, the map $F_\veps$ is a contraction on each one of
\begin{equation}
\Uplambda^m_\veps:=\big\{\widetilde{\Lambda}_\veps\in \widetilde{W}^{1}_\infty([0,(m+1)T]):\,0\!\le\!\widetilde{\Lambda}_\veps'\!\le\! \|f_\veps\|_\infty/\veps\;\text{and}\;\widetilde{\Lambda}_\veps|_{[0,mT]}\!=\!\Lambda_\veps|_{[0,mT]}\big\},\;m\ge 0. 
\end{equation}
\end{prop}

\noindent\textbf{Proof. Step 1.} We claim that, for any $\Lambda_\veps\in \widetilde{W}^{1}_\infty([0,T])$ with $\|\Lambda_\veps'\|_\infty\le\|f_\veps\|_\infty/\veps$, the function $F_\veps(\Lambda_\veps)$ is non-decreasing with a Lipschitz constant of at most  $\|f_\veps\|_\infty/\veps$. Indeed, let $p_\veps\in W^{1,2}_2([0,T]\times[0,\infty))$ be the solution of \eqref{ppde}. By the representation in Proposition \ref{prop:FK} we easily obtain the lower bound $p_\veps\ge 0$. For an upper bound, we consider $(t,x)\in(0,T]\times(0,\infty)$ and repeat \cite[proof of Lemma 3.1, Step 2 until (3.9)]{NS} literally for $p_\veps(t-s,y)$ instead of their $\zeta(s,y)$ to infer
\begin{equation}
p_\veps(t,x)=\ev[p_\veps(t-\overline{\tau},\overline{Z}_{\overline{\tau}})]=\ev\big[f_\veps(\overline{Z}_t)\,\mathbf{1}_{\{\overline{\tau}=t\}}\big]
+\ev\big[p_\veps(t-\overline{\tau},0)\,\mathbf{1}_{\{\overline{\tau}<t\}}\big],
\end{equation}
where $\overline{Z}_s:=x+B_s-\Lambda_\veps(t\!-\!s)+\Lambda_\veps(t)$, $s\in[0,t]$ and $\overline{\tau}:=\inf\{s\in[0,t]:\overline{Z}_s\!=\!0\}\!\wedge t$. So,
\begin{equation}\label{max_prin}
p_\veps(t,x)\le \|f_\veps\|_\infty\vee \max_{s\in[0,t]} p_\veps(s,0),\;(t,x)\in(0,T]\times(0,\infty).
\end{equation}
But $\max_{s\in[0,t]} p_\veps(s,0)>\|f_\veps\|_\infty$ yields $\partial_x p_\veps(s^*,0)>0$ for $s^*\in\arg\max_{s\in[0,t]} p_\veps(s,0)\!\neq\!\varnothing$ by the Robin boundary condition in \eqref{ppde} and $2\Lambda_\veps'(s^*)\le 2\|f_\veps\|_\infty/\veps<\frac{\alpha}{\veps}$, and we deduce
\begin{equation}
\sup_{(0,t]\times(0,\infty)} p_\veps>p_\veps(s^*,0)=\max_{s\in[0,t]} p_\veps(s,0)=\|f_\veps\|_\infty\vee\max_{s\in[0,t]} p_\veps(s,0),
\end{equation}
a direct contradiction to \eqref{max_prin}. Thus, we end up with $0 \le p_\veps\le \|f_\veps\|_\infty$. In conjunction with \eqref{Fdef}, this results in
\begin{equation}
0\le  F_\veps(\Lambda_\veps)(t_2)-F_\veps(\Lambda_\veps)(t_1) \le \frac{\|f_\veps\|_\infty}{\veps}\,(t_2-t_1),\;0\le t_1<t_2\le T,
\end{equation}
establishing the claim.

\medskip

\noindent\textbf{Step 2.}~We aim to apply the Banach fixed point theorem in the complete metric space $\Uplambda^0_\veps$. Take $\Lambda_{\veps,1},\Lambda_{\veps,2}\in\Uplambda^0_\veps$ and the associated solutions $p_{\veps,1},p_{\veps,2}\in W^{1,2}_2([0,T]\times[0,\infty))$ of \eqref{ppde}. Given $t\in(0,T]$, let $\eta_\delta\in W^{1,2}_2([0,t]\times[0,\infty))$ solve
\begin{align}\label{eta_pde}
\begin{cases}
\;\pp_s\eta_\delta+\frac{1}{2}\pp_{xx}\eta_\delta-\Lambda'_{\veps,1}\,\pp_x\eta_\delta=0\;\;\text{on}\;\;[0,t]\times[0,\infty), \\
\;\eta_\delta(t,\cdot)=\xi_\delta\;\;\text{and}\;\;\alpha\eta_\delta(\cdot,0)=\veps\pp_x\eta_\delta(\cdot,0),
\end{cases}
\end{align}
where $\xi_\delta(\cdot)=\xi(\cdot/\delta)/\delta$, $\delta>0$ for some non-negative $\xi\in C^\infty_c((0,\infty))$ that integrates to one. Arguing as in the derivation of \eqref{duality} we find
\begin{equation}\label{pdelta}
\begin{split}
&\,\left|\int_{0}^\infty \big(p_{\veps,1}(t,x)-p_{\veps,2}(t,x)\big)\,\xi_\delta(x)\,\dd x\right|\\
&=\left|\int_0^t \int_{0}^\infty \left(\Lambda'_{\veps,2}(s)-\Lambda'_{\veps,1}(s)\right)\pp_x\eta_\delta(s,x)\,p_{\veps,2}(s,x)\,\dd x\,\dd s\right|\\
&\le \|f_\veps\|_\infty\,\sup_{s\in[0,t]}\left|\Lambda'_{\veps,2}(s)-\Lambda'_{\veps,1}(s)\right|\int_0^t \int_{0}^\infty |\pp_x\eta_\delta(s,x)|\,\dd x\,\dd s.
\end{split}
\end{equation}

\smallskip

To bound $\int_0^t \int_{0}^\infty |\pp_x\eta_\delta(s,x)|\,\dd x\,\dd s$ we use \cite[Chapter VI, Theorem 1.10]{garroni1992green}  and represent $\eta_\delta$ by means of a Green's function $G$:
\begin{equation}
\eta_\delta(s,x)=-\int_s^t\int_0^\infty G(s,x;r,y)\,\Lambda'_{\veps,1}(r)\,\partial_x\eta_\delta(r,y)\,\mathrm{d}y\,\mathrm{d}r
+\int_0^\infty G(s,x;t,y)\,\xi_\delta(y)\,\dd y.
\end{equation}
Differentiating in $x$ and employing the estimates $0\le\Lambda'_{\veps,1}\le\|f_\veps\|_\infty/\veps$ and
\begin{align}
|\pp_x G(s,x;r,y)| \le C(r-s)^{-1}\,e^{-c(y-x)^2/(r-s)},
\end{align}
with $C=C(\veps,\alpha)<\infty$, $c=c(\veps,\alpha)>0$, (see \cite[Chapter VI, Theorem 1.10(i)]{garroni1992green}) we get   
\begin{equation}
\begin{split}
|\partial_x \eta_\delta(s,x)|\le &\;\frac{\|f_\veps\|_\infty}{\veps}\int_s^t \int_0^\infty  C(r-s)^{-1}\,e^{-c(y-x)^2/(r-s)}\,|\partial_x \eta_\delta(r,y)|\,\mathrm{d}y\,\mathrm{d}r \\
&\;+\int_0^\infty C(t-s)^{-1}\,e^{-c(y-x)^2/(t-s)}\,\xi_\delta(y)\,\mathrm{d}y.
\end{split}
\end{equation}
Next, we integrate in $(s,x)$ and control the $\mathrm{d}x$ integrals over $[0,\infty)$ by those over $\R$: 
\begin{equation}
\begin{split}
\int_0^t \int_0^\infty |\partial_x \eta_\delta(s,x)|\,\mathrm{d}x\,\mathrm{d}s
\le &\;\frac{\|f_\veps\|_\infty}{\veps}\,\frac{2C\sqrt{\pi}}{\sqrt{c}}\int_0^t\int_0^\infty \sqrt{r}\,|\partial_x \eta_\delta(r,y)|\,\mathrm{d}y\,\mathrm{d}r \\
&\;+\frac{2C\sqrt{\pi}}{\sqrt{c}}\int_0^\infty \sqrt{t}\,\xi_\delta(y)\,\mathrm{d}y.
\end{split}
\end{equation}
For $0<T<\frac{\veps^2}{\|f_\veps\|_\infty^2}\,\frac{c}{4C^2\pi}$, it follows via $\int_0^\infty \xi_\delta(y)\,\mathrm{d}y=1$ that
\begin{equation}\label{etax}
\int_0^t \int_0^\infty |\partial_x \eta_\delta(s,x)|\,\mathrm{d}x\,\mathrm{d}s\le\frac{\frac{2C\sqrt{\pi }}{\sqrt{c}}\,\sqrt{T}}{1-\frac{\|f_\veps\|_\infty}{\veps}\,\frac{2C\sqrt{\pi }}{\sqrt{c}}\,\sqrt{T}}. 
\end{equation}

\smallskip

Recalling the definition of $F_\veps$ in \eqref{Fdef}, inserting \eqref{etax} into \eqref{pdelta}, and taking the limit $\delta\downarrow0$ and then the supremum over $t\in[0,T]$ we arrive at 
\begin{equation}
\begin{split}
\|F_\veps(\Lambda_{\veps,1})'-F_\veps(\Lambda_{\veps,2})'\|_\infty
=&\;\frac{1}{\veps} \sup_{t\in[0,T]} |p_{\veps,1}(t,0)-p_{\veps,2}(t,0)| \\
\le  &\;\frac{\|f_\veps\|_\infty}{\veps}\,\|\Lambda'_{\veps,2}-\Lambda'_{\veps,1}\|_\infty\,\frac{\frac{2C\sqrt{\pi }}{\sqrt{c}}\,\sqrt{T}}{1-\frac{\|f_\veps\|_\infty}{\veps}\,\frac{2C\sqrt{\pi }}{\sqrt{c}}\,\sqrt{T}}
\end{split}
\end{equation}
for $0<T<\frac{\veps^2}{\|f_\veps\|_\infty^2}\,\frac{c}{4C^2\pi}$. Thus, for all small enough $T=T(\|f_\veps\|_\infty,\veps,\alpha)>0$, the map $F_\veps$ is a contraction on $\Uplambda^0_\veps$ and possesses a unique fixed point therein.  

\medskip

\noindent\textbf{Step 3.} We conclude by using induction over $m\ge0$ to show that, for all small enough $T=T(\|f_\veps\|_\infty,\veps,\alpha)>0$ as in Step 2, the map $F_\veps$ is a contraction on each one of $\Uplambda^m_\veps$, $m\ge0$ and the resulting fixed points in $\Uplambda^m_\veps$, $m\ge0$ give the respective unique fixed points $\Lambda_\veps$ of $F_\veps$ in $\widetilde{W}^{1}_\infty([0,(m+1)T])$, $m\ge0$ with $\|\Lambda'_\veps\|_\infty\le\|f_\veps\|_\infty/\veps$. Since we have already established the statement for $m=0$, we turn to the induction step for $m\ge1$. Pick $\Lambda_{\veps,1},\Lambda_{\veps,2}\in\Uplambda^m_\veps$, the associated solutions $p_{\veps,1},p_{\veps,2}\in W^{1,2}_2([0,(m+1)T]\times[0,\infty))$ of \eqref{ppde}, $t\in(mT,(m+1)T]$, and $\eta_\delta\in W^{1,2}_2([0,t]\times[0,\infty))$ solving \eqref{eta_pde}. Repeating \eqref{pdelta} and exploiting $\Lambda_{\veps,1}|_{[0,mT]}=\Lambda_{\veps,2}|_{[0,mT]}$ we deduce  
\begin{equation}
\begin{split}
& \left|\int_{0}^\infty \big(p_{\veps,1}(t,x)-p_{\veps,2}(t,x)\big)\,\xi_\delta(x)\,\dd x\right| \\
& \le \|f_\veps\|_\infty\,\sup_{s\in[mT,t]} \left|\Lambda'_{\veps,2}(s)-\Lambda'_{\veps,1}(s)\right|\int_{mT}^t \int_{0}^\infty |\pp_x\eta_\delta(s,x)|\,\dd x\,\dd s.
\end{split}
\end{equation}
As before, the latter double integral cannot exceed the right-hand side of \eqref{etax}. Taking the limit $\delta\downarrow0$ and then the supremum over $t\in(mT,(m+1)T]$ we get via $F_\veps(\Lambda_{\veps,1})|_{[0,mT]}=F_\veps(\Lambda_{\veps,2})|_{[0,mT]}$:
\begin{equation}
\begin{split}
\|F_\veps(\Lambda_{\veps,1})'-F_\veps(\Lambda_{\veps,2})'\|_\infty
=&\;\frac{1}{\veps}\sup_{t\in(mT,(m+1)T]} |p_{\veps,1}(t,0)-p_{\veps,2}(t,0)| \\
\le&\;   \frac{\|f_\veps\|_\infty}{\veps}\,\|\Lambda'_{\veps,2}-\Lambda'_{\veps,1}\|_\infty\,\frac{\frac{2C\sqrt{\pi }}{\sqrt{c}}\,\sqrt{T}}{1-\frac{\|f_\veps\|_\infty}{\veps}\,\frac{2C\sqrt{\pi }}{\sqrt{c}}\,\sqrt{T}}. 
\end{split} 
\end{equation} 
The contraction property of $F_\veps$ on $\Uplambda^m_\veps$ readily follows and upon combining it with the induction hypothesis we complete the induction step. \qed

%%%%%%%%%%%%%%%%%%%%
\section{Monotonicity in $\veps$}\label{sec:mono}
%%%%%%%%%%%%%%%%%%%%

The next proposition asserts that the free boundaries $\{\Lambda_\veps\}_{\veps>0}$, obtained via Proposition \ref{PDEsol}, increase pointwise as $\veps\downarrow0$.   

\begin{prop}\label{prop:mono}
Let $f$ be non-negative and bounded, with $\int_0^\infty f(x)\,\mathrm{d}x=1$ and $\|f\|_\infty<\frac{\alpha}{2}$. Then, for any $\veps>\widehat{\veps}>0$, it holds $\Lambda_\veps(t)\le\Lambda_{\widehat{\veps}}(t)$ for all $t\ge0$, where $\Lambda_\veps$ is the free boundary in the problem \eqref{stefan_reg} with $\|\Lambda'_\veps\|_\infty\le\|f_\veps\|_\infty/\veps$ and 
$\Lambda_{\widehat{\veps}}$ is the respective free boundary for $\widehat{\veps}$.~In particular, $\{\Lambda_\veps\}_{\veps>0}$ tend pointwise to some $\Lambda_0$ as $\veps\downarrow0$. 
\end{prop}

\noindent\textbf{Proof.} We fix $\veps>\widehat{\veps}>0$ and take a small enough $T>0$ as in the second statement of Proposition \ref{PDEsol}, so that the maps $F_\veps$ and $F_{\widehat{\veps}}$ are contractions on each one of the spaces $\Uplambda^m_\veps$, $m\ge0$ and $\Uplambda^m_{\widehat{\veps}}$, $m\ge0$, respectively, defined therein. We argue by induction over $m\ge0$ that $\Lambda_\veps(t)\le\Lambda_{\widehat{\veps}}(t)$ for all $t\in[0,(m+1)T]$. For $m=0$, we use the contraction properties of $F_\veps$ and $F_{\widehat{\veps}}$ on $\Uplambda^0_\veps$ and $\Uplambda^0_{\widehat{\veps}}$, respectively, to conclude that $\Lambda_\veps=\lim_{n\to\infty} \Lambda^n_\veps$ in $\Uplambda^0_\veps$ and $\Lambda_{\widehat{\veps}}=\lim_{n\to\infty} \Lambda^n_{\widehat{\veps}}$ in $\Uplambda^0_{\widehat{\veps}}$, where 
\begin{equation}
\Lambda^0_\veps\equiv 0,\;\;\Lambda^n_\veps=F_\veps(\Lambda^{n-1}_\veps),\;n\ge 1\quad\text{and}\quad
\Lambda^0_{\widehat{\veps}}\equiv 0,\;\;\Lambda^n_{\widehat{\veps}}=F_{\widehat{\veps}}(\Lambda^{n-1}_{\widehat{\veps}}),\;n\ge 1. 
\end{equation}
In view of \eqref{eq:sformula} and Proposition \ref{prop:burdzy}, we can employ induction over $n\ge0$ to establish $\Lambda^n_\veps(t)\le\Lambda^n_{\widehat{\veps}}(t)$, $t\in[0,T]$ for all $n\ge0$. Consequently, we have $\Lambda_\veps(t)\le\Lambda_{\widehat{\veps}}(t)$, $t\in[0,T]$. 

\medskip

For the induction step, we let $m\ge1$ and rely on the contraction properties of $F_\veps$ and $F_{\widehat{\veps}}$ on $\Uplambda^m_\veps$ and $\Uplambda^m_{\widehat{\veps}}$, respectively, to deduce that $\Lambda_\veps=\lim_{n\to\infty} \Lambda^n_\veps$ in $\Uplambda^m_\veps$ and $\Lambda_{\widehat{\veps}}=\lim_{n\to\infty} \Lambda^n_{\widehat{\veps}}$ in $\Uplambda^m_{\widehat{\veps}}$, where 
\begin{equation}
\begin{split}
&\Lambda^0_\veps(t)=\begin{cases}
\;\Lambda_\veps(t) &\text{if}\quad t\in[0,mT], \\
\;\Lambda_\veps(mT) &\text{if}\quad t\in(mT,(m+1)T],	 
\end{cases}
\quad\Lambda^n_\veps=F_\veps(\Lambda^{n-1}_\veps),\;n\ge 1\quad\text{and}\quad \\
&\Lambda^0_{\widehat{\veps}}(t)=\begin{cases}
\;\Lambda_{\widehat{\veps}}(t) &\text{if}\quad t\in[0,mT], \\
\;\Lambda_{\widehat{\veps}}(mT)& \text{if}\quad t\in(mT,(m+1)T],	 
\end{cases}
\quad\Lambda^n_{\widehat{\veps}}=F_{\widehat{\veps}}(\Lambda^{n-1}_{\widehat{\veps}}),\;n\ge 1. 
\end{split}
\end{equation}
Putting the induction hypothesis together with \eqref{eq:sformula} and Proposition \ref{prop:burdzy} we infer $\Lambda^n_\veps(t)\le\Lambda^n_{\widehat{\veps}}(t)$, $t\in[0,(m+1)T]$ inductively over $n\ge0$. Thus, $\Lambda_\veps(t)\le\Lambda_{\widehat{\veps}}(t)$ for all $t\in[0,(m+1)T]$. Finally, the pointwise convergence of $\{\Lambda_\veps\}_{\veps>0}$ as $\veps\downarrow0$ follows from this and the uniform boundedness of $\{\Lambda_\veps\}_{\veps>0}$ (which is immediate from \eqref{eq:sformula}). \qed

%%%%%%%%%%%%%%%%%%%%
\section{Identification of the limit}
%%%%%%%%%%%%%%%%%%%%

We define the right-continuous modification $\widetilde{\Lambda}(t)=\lim_{s\downarrow t} \Lambda_0(s)$, $t\ge0$ of the pointwise limit $\Lambda_0$ from Proposition \ref{prop:mono} and show that $\widetilde{\Lambda}$ solves the limiting problem \eqref{limprob}.

\begin{prop}\label{prop:limprob}
Let $f$ be non-negative and bounded, with $\int_0^\infty f(x)\,\mathrm{d}x=1$ and $\|f\|_\infty<\frac{\alpha}{2}$. Then, the right-continuous modification $\widetilde{\Lambda}(t)=\lim_{s\downarrow t} \Lambda_0(s)$, $t\ge0$ of the pointwise limit $\Lambda_0$ from Proposition \ref{prop:mono} is a solution of the problem \eqref{limprob}.
\end{prop}

\noindent\textbf{Proof.} With $X_t:=X_{0-}+B_t-\widetilde{\Lambda}(t)$, $t\ge0$ and $\tau:=\inf\{t\ge0:\,X_t\le0\}$, we need to check that $\alpha\widetilde{\Lambda}(t)=2\mathbb{P}(\tau\le t)$, $t\ge0$. Since $\widetilde{\Lambda}$, $t\mapsto\mathbb{P}(\tau\le t)$ are right-continuous and the set of continuity points of $\widetilde{\Lambda}$ is dense in $[0,\infty)$, we restrict our attention to the continuity points $t$ of $\widetilde{\Lambda}$ throughout and note that $\widetilde{\Lambda}(t)=\Lambda_0(t)$ for such $t$. Moreover, by Proposition \ref{prop:mono} we have $\Lambda_\veps\le \Lambda_0\le \widetilde{\Lambda}$ for all $\veps>0$, and thus $\{\tau>t\}$ implies $\min_{0\le s\le t} (X^\veps_0+B_s-\Lambda_\veps(s))>0$ and $L^{\veps,0}_t=0$ for all $\veps>0$. Consequently, for all $\veps>0$,
\begin{equation}
\Lambda_\veps (t)
=\frac{2}{\alpha} \big(1-\mathbb{E}[\mathbf{1}_{\{\tau> t\}}]-\mathbb{E}\big[\mathbf{1}_{\{\tau\leq t\}}\,e^{-\alpha L_t^{\veps,0}/\veps}\big]\big)=\frac{2}{\alpha} \big(\mathbb{P}(\tau\le t)-\mathbb{E}\big[\mathbf{1}_{\{\tau\leq t\}}\,e^{-\alpha L_t^{\veps,0}/\veps}\big]\big).
\end{equation}

\smallskip

It remains to show that $\lim_{\veps\downarrow0} \mathbb{E}\big[\mathbf{1}_{\{\tau\leq t\}}\,e^{-\alpha L_t^{\veps,0}/ \veps}\big]=0$. Let $L_t:=\lim_{\veps\downarrow0}L_t^{\veps,0}\in[0,\infty]$, $t\ge0$, which exists thanks to Proposition \ref{prop:burdzy}. Recall that \eqref{SkorFormula} and Proposition \ref{prop:refl} give
\begin{align}
L_t^{\veps,0} = \Big(\!-\!\min_{0\le s\le t} \big(X^\veps_0+B_s-\Lambda_\veps(s)\big)\!\Big)\vee 0. 
\end{align}
It follows that 
\begin{align}\label{eq:lambdaskorokhod}
L_t=\lim_{\veps\downarrow0} \Big(\!-\!\min_{0\le s\le t} \big(X^\veps_0+B_s-\Lambda_\veps(s)\big)\!\Big)\vee 0
=\Big(\!-\!\min_{0\le s\le t} \big(X_{0-}+B_s-\widetilde{\Lambda}(s)\big)\!\Big)\vee 0,
\end{align}
where $\min_{0\le s\le t} (X_{0-}+B_s-\widetilde{\Lambda}(s))$ is well-defined and a.s.~attained at a continuity point of $\widetilde{\Lambda}$ since $s\mapsto X_{0-}+B_s-\widetilde{\Lambda}(s)$ is lower semi-continuous, $t$ is a continuity point of $\widetilde{\Lambda}$, standard Brownian motion a.s.~instantaneously enters into the negative half-line (cf.~\cite[Chapter 2, Theorem 9.23(ii)]{karatzas1998brownian}), and $\widetilde{\Lambda}$ is non-decreasing. On $\{\tau\le t\}$, it holds $\min_{0\le s\le t} (X_{0-}+B_s-\widetilde{\Lambda}(s))<0$ a.s., hence also $L_t>0$ a.s., due to $\{\tau=t\}=\{X_{0-}+B_t-\widetilde{\Lambda}(t)=0\}$ being a $\bP$-null set, the fact that standard Brownian motion a.s.~instantaneously enters into the negative half-line, and the monotonicity of $\widetilde{\Lambda}$. Thus, $\lim_{\veps\downarrow0} \mathbf{1}_{\{\tau\leq t\}}\,e^{-\alpha L_t^{\veps,0}/ \veps}=0$ a.s., yielding the desired $\lim_{\veps\downarrow0} \mathbb{E}\big[\mathbf{1}_{\{\tau\leq t\}}\,e^{-\alpha L_t^{\veps,0}/ \veps}\big]=0$ via the dominated convergence theorem. \qed
	
\medskip	

To complete the proof of Theorem \ref{thm:main} it remains to identify the pointwise limit $\Lambda_0$ from Proposition \ref{prop:mono} as the unique solution $\Lambda$ of \eqref{limprob} from Proposition \ref{thm:unique}. 

\begin{prop}\label{prop:unique}
Let $f$ be non-negative and bounded, with $\int_0^\infty f(x)\,\mathrm{d}x=1$ and $\|f\|_\infty<\frac{\alpha}{2}$. Then, the pointwise limit $\Lambda_0$ from Proposition \ref{prop:mono} agrees with the unique solution $\Lambda$ of \eqref{limprob} from Proposition \ref{thm:unique}.
\end{prop}

\noindent\textbf{Proof.} Combining Propositions \ref{prop:limprob} and \ref{thm:unique} we deduce  $\widetilde{\Lambda}=\Lambda$ and the continuity of $\widetilde{\Lambda}$. Consequently, $\Lambda_0=\widetilde{\Lambda}=\Lambda$. \qed

%\appendix
%\section*{Appendix}

\bigskip\bigskip

\bibliographystyle{alpha}
\bibliography{biblio}

\bigskip\bigskip

\end{document}